\documentclass{siamltex}

\usepackage{amsfonts,latexsym,amsmath,graphicx,psfrag}
\usepackage{amssymb}

\newtheorem{remark}[theorem]{Remark}

\newcommand{\notmid}{\mid\kern-0.5em\not\kern0.5em}

\newcommand{\re}[1]{(\ref{#1})}
\newcommand{\wt}[1]{\widetilde{#1}}
\newcommand{\wh}[1]{\widehat{#1}}

\newcommand{\pdpd}[2]{\frac{\partial #1}{\partial #2}}

\newcommand{\dd}{\mathrm{d}}

\newcommand{\ignore}[1]{}

\newcommand{\mstrut}[1]{\mbox{\rule{0mm}{#1}}}

\newcommand{\barr}{\begin{array}}
\newcommand{\earr}{\end{array}}

\def\bfo{\begin{eqnarray*}}
\def\efo{\end{eqnarray*}}
\def\ba{\begin{eqnarray*}}
\def\ea{\end{eqnarray*}}
\def\beq{\begin{eqnarray}}
\def\eeq{\end{eqnarray}}

\def\hat{\widehat}
\def\tilde{\widetilde}

\def\T{{\mathcal T}}

\def\bra{\langle}
\def\cet{\rangle}

\def\p{\partial}

\title{Reconstruction of a conformally Euclidean metric from local boundary diffraction travel times\thanks{This research was supported by National Science Foundation grant CMG DMS-1025318, the members of the Geo-Mathematical Imaging Group at Purdue University, the Finnish Centre of Excellence in Inverse Problems Research, Academy of Finland project COE 250215, the Research Council of Norway, and the VISTA project. The research was initialized at the Program on Inverse Problems and Applications at MSRI, Berkeley, during the Fall 2010.}}
  
  \author{Maarten V. de Hoop\thanks{Department of Mathematics, Purdue University 150 N. University Street, West Lafayette IN 47907, USA ({\tt mdehoop@math.purdue.edu})}
\and  
Sean F. Holman\thanks{Department of Mathematics, Purdue University 150 N. University Street, West Lafayette IN 47907, USA ({\tt sfholman@math.purdue.edu}), corresponding author.}
 \and
Einar Iversen\thanks{NORSAR, Gunnar Randers vei 15, P.O. Box 53, 2027 Kjeller, Norway ({\tt Einar.Iversen@norsar.com})}
\and
Matti Lassas\thanks{Department of Mathematics and Statistics, Gustaf Hallstromin katu 2b, FI-00014 University of Helsinki, Helsinki, Finland ({\tt matti.lassas@helsinki.fi})}
\and
Bj{\o}rn Ursin\thanks{Department of Petroleum Engineering and Applied Geophysics, Norwegian University of Science and Technology, S.P. Andersensvei 15A, NO-7491 Trondheim, Norway ({\tt bjorn.ursin@ntnu.no})}
}

\begin{document}

\maketitle

\begin{abstract}
We consider a region $M$ in $\mathbb{R}^n$ with
boundary $\partial M$ and a metric $g$ on $M$ conformal to the Euclidean
metric. We analyze the inverse problem, originally formulated by Dix~\cite{dix}, of reconstructing $g$ from boundary measurements
associated with the single scattering of seismic waves in this region. In our formulation the measurements determine the shape operator of wavefronts outside of $M$ originating at diffraction points within $M$. We develop an explicit
reconstruction procedure which consists of two steps. In the first
step we reconstruct the directional curvatures and the metric in what are essentially Riemmanian normal coordinates; in the second step we develop a conversion to Cartesian
coordinates. We admit the presence of conjugate points.  In dimension $n \geq 3$ both steps involve the solution of a system of ordinary differential equations. In dimension $n=2$ the same is true for the first step, but the second step requires the solution of a Cauchy problem for an elliptic operator which is unstable in general. The first step of the procedure applies for general metrics.
\end{abstract}



 
\pagestyle{myheadings}
\thispagestyle{plain}
\markboth{
        }{Reconstruction of a Riemannian metric}

\section{Introduction}

We consider a region, $M$, in $\mathbb{R}^n$ with a smooth boundary $\partial M$. We assume that there is a Riemannian metric, $g$, on $M$ that is conformal to the Euclidean metric with conformal factor $v^{-2}$ where $v \in C^\infty(\overline{M})$ is strictly positive. This means that $g(x) = v^{-2} \mathbf{e}$, where $\mathbf{e}$ is the Euclidean metric, or, in Cartesian coordinates $x = (x^1,\ldots,x^n)$, $g^{ij}(x) = v(x)^2 \delta^{ij}$. We analyze the inverse problem of reconstructing $g$ based on measurements of the curvature of wavefronts produced by point diffractors located inside $M$ and propagated according to the wave operator on $(M,g)$. Indeed, geodesics for the metric $g$ are rays following the propagation of singularities by a parametrix corresponding to this wave operator. In the seismic context $v$ is the wave speed. As originally formulated by Dix \cite{dix}, this type of data may in some cases be reconstructed from reflection data by variation of source and receiver locations at the surface of the earth. In particular, it is possible to recover from reflection data the shape operator for the wave front produced by a given point diffractor where the rays beginning at the diffractor intersect $\partial M$ orthogonally \cite{GJI}. Dix developed a procedure, with a formula, for reconstructing one-dimensional wave speed profiles in a half space from reflection data. Since Dix various adaptations have been considered to admit more general wave speed functions in a half space. We mention the work of Shah \cite{shah}, Hubral \& Krey \cite{hubralK}, Dubose, Jr. \cite{dubose}, and Mann \cite{mann}. We consider here the case of higher dimensional regions with Riemannian metrics conformal to the Euclidean metric. Our method is different in the cases of $n=2$ and $n\geq 3$ and in fact we expect better results in the case $n\geq 3$. The problem is closely related to the problem where broken geodesics are observed on the boundary \cite{KKL} or when the Cauchy data of the solution of the wave equations are observed on the boundary, see e.g. \cite{AKKLT,BK,KKL} and references therein.

Assuming that we know $v$ and all of its derivatives on $\partial M$ (c.f. \cite{LassasSU:2003}), we may extend $v$ to a function, which we will also denote by $v$, on a complete manifold $\widetilde{M}$ compactly containing $M$ (in fact we can take $\widetilde{M} = \mathbb{R}^n$). The corresponding extended metric is also denoted simply be $g$. As described in detail below, we measure the curvature of generalized spheres for $g$ centered at ``diffraction'' points intersected with an open subset of $\widetilde{M} \setminus M$. From these data, and assuming we know $v$ in $\widetilde M \setminus M$, we show an explicit method to determine the function $v$ in the Cartesian coordinates $x = (x^1,\ldots,x^n)$ along geodesics of $g$ which connect the diffractions points to the measurement region. This method can be viewed as a generalization of the work of Iversen and Tygel \cite{tygelI}. We now proceed to introduce the concepts and notations necessary for the statement of our main results.

For any $(x,\eta) \in \Omega \widetilde{M}$ ($\Omega$ indicates the unit sphere bundle with respect to $g$) we will write $\gamma_{x,\eta}$ for the geodesic with initial data $\gamma_{x,\eta}(0) = x$, $\dot \gamma_{x,\eta}(0) = \eta$. For $r$ in the domain of $\gamma_{x,\eta}$ we let $\mathcal{C}_r(x,\eta)$ denote the set of times $t$ such that $t$ is conjugate to $r$ along $\gamma_{x,\eta}$ (by this we mean that $(t-r) \dot \gamma_{x,\eta}(r)$ is a critical point for $\mathrm{exp}_{\gamma_{x,\eta}(r)}$). 

For the moment we fix $(x_0,\eta_0) \in \Omega (\widetilde{M} \setminus M)$ and use the notation $\mathcal{C}_r$ for $\mathcal{C}_r(x_0,\eta_0)$. We will refer to the image of the set
\[
\{\xi \in T_{y} \tilde M \ : \ |\xi|_g = R\}
\]
under the exponential map as the {\it generalized sphere of radius $R$ centered at $y$}. When $t>r \geq 0$ is in the domain of $\gamma_{x_0,\eta_0}$ and $t \notin \mathcal{C}_r$
there is a small portion of the generalized sphere of radius $t-r$ centered at $y_t :=\gamma_{x_0,\eta_0}(t)$ containing
$\gamma_{x_0,\eta_0}(r)$ that is an embedded submanifold $\Sigma_{r,t}$
of $\widetilde{M}$. Indeed, in this case we can define a vector field
$\nu_{r,t}$ in a neighborhood of $\gamma_{x_0,\eta_0}(r)$ by writing
for $\xi \in T_{y_t} \widetilde{M}$ in a small neighborhood of
$-(t-r)\dot \gamma_{x_0,\eta_0}(t)$:
\[
   \nu_{r,t}(\mathrm{exp}_{y_t}(\xi))
   = \frac{1}{|\xi|_g} \left. \frac{\partial}{\partial s}
          \right |_{s=1} \mathrm{exp}_{y_t}(s\xi) .
\]
Geometrically, $\nu_{r,t}$ gives the outward pointing normal vector
fields to a part of the generalized sphere centered at $y_t$ near
$\gamma_{x_0,\eta_0}(r)$. The shape operator $S_{r,t} \in
(T_1^1)_{\gamma_{x_0,\eta_0}(r)} \widetilde{M}$ of $\Sigma_{r,t}$ at $\gamma_{x_0,\eta_0}(r)$ is
then given by
\[
   S_{r,t} X = \nabla_X \nu_{r,t}
\]
for all $X \in T_{\gamma_{x_0,\eta_0}(r)} \tilde{M}$, where $\nabla$
is the Levi-Civita connection for $g$. For the reconstruction of $v$ we assume
that $S_{0,t}$ is known for all $t>0$ such that $t \notin
\mathcal{C}_0$. In reflection seismology one refers to $\Sigma_{r,t}$
as the (partial) front of a point diffractor located at $y_{t}$.

\begin{figure}\label{notFig}
\center{
\includegraphics{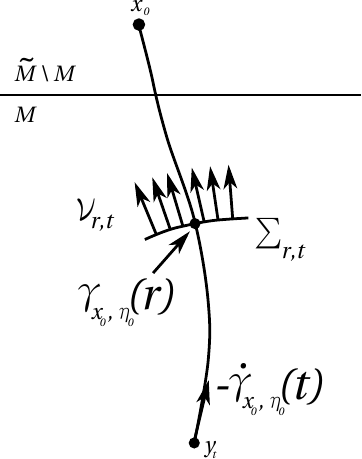}}
\caption{An illustration of the various notations introduced.}
\end{figure}

\begin{figure}\label{fig:coord}
\center{\includegraphics[scale=.5]{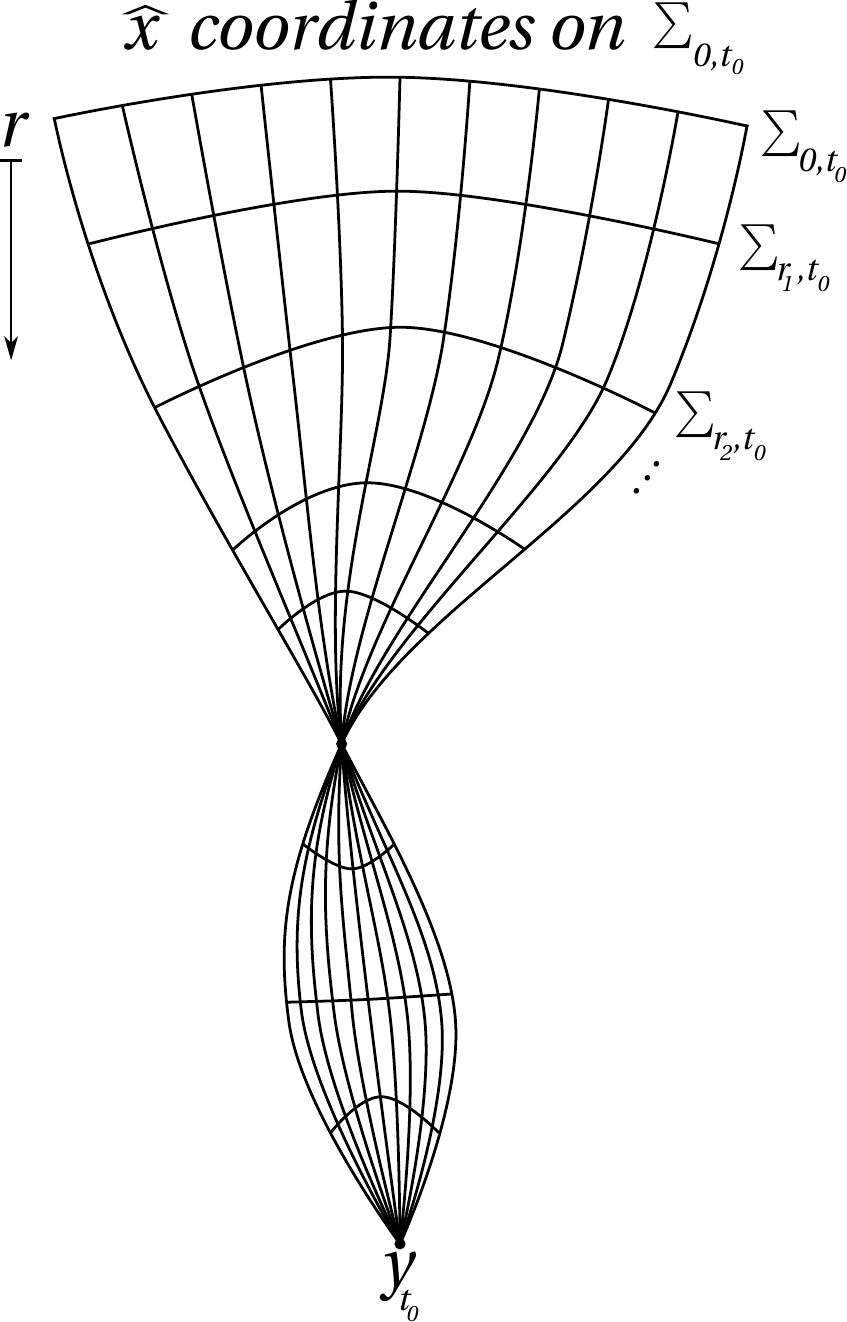}}
\caption{A figure illustrating the coordinates $(\hat{x},r)$ defined on the set $W(t_0)$.}
 \label{WFig}
\end{figure}

We now introduce a mapping which defines local coordinates in which we will perform our
initial calculations. We begin with picking a large $t_0 >0$ in the domain of $\gamma_{x_0,\eta_0}$ such that
$t_0 \notin \mathcal{C}_0$. Next, let us take local coordinates
$\hat{x} = (\hat{x}^1,\ldots,\hat{x}^{n-1})$ on $\Sigma_{0,t_0}$ such
that $(\hat{x}^1,\ldots,\hat{x}^{n-1})=0$ defines $x_0$ and suppose $\Phi_{t_0}: \Sigma_{0,t_0} \mapsto U\subset \mathbb{R}^{n-1}$ is the coordinate map. We will assume that the image of these maps, $U$, is the same for all $t_0$. For $\hat{x} \in U$, let $\gamma_{\hat{x}}^{t_0}$ be the
geodesic with a special choice of initial data: $\dot\gamma^{t_0}_{\hat{x}}(0) = -
\nu_{0,t_0}(\Phi_{t_0}^{-1}(\hat{x}))$. Then we define coordinates, $(\hat{x},r)$, on
some set by the inverse of the map 
\begin{equation}\label{map}
  \Psi_{t_0}(\hat{x},r) = \gamma^{t_0}_{\hat{x}}(r).
\end{equation}
We define
\[
W(t_0) = \Psi_{t_0}\left ( \left \{ (\hat{x},r) \in U \times \mathbb{R} \ : \ r < t_0, \ r \notin \mathcal{C}_{t_0}\left ( \Phi^{-1}_{t_0}(\hat{x}), -\nu_{0,t_0}(\Phi^{-1}_{t_0}(\hat{x})) \right )  \right \} \right )
\]
and so that the map $\Psi_{t_0}$ defines a local parametrization on $W(t_0)$. It may not be a global parametrization because it may not be injective. The $(\hat{x},r)$ local
coordinates on $W(t_0)$, basically, are Riemannian normal coordinates centered at
$y_{t_0}$, but parametrized in a particular way: $\hat{x}$ can be
thought of as a parametrization of part of the sphere of
radius $t_0$ in $T_{y_{t_0}} \tilde{M}$, and then $r$ corresponds to
the radial variable in $T_{y_{t_0}} \tilde{M}$. Figure \ref{WFig} shows a possible depiction of $W(t_0)$ and the coordinates defined there. We note that the
domain $W(t_0)$ includes $\gamma_{x_0,\eta_0}([0,t_0))
\setminus \{\gamma_{x_0,\eta_0}(r) \, : \, r \in \mathcal{C}_{t_0}
\}$.
We also define for some $L >0$
\begin{equation}\label{Wdef}
W = \bigcup_{t_0 \in (0,L)} W(t_0).
\end{equation}
In our reconstruction we cover $M$ by sets like $W(t_0)$, and then recover $v$ on each of these regions separately. In the case that there are conjugate points to $t_0$ along $\gamma_{x_0,\eta_0}$ we will also have to cover some regions with more than one such set. 
Note that along $\gamma_{x_0,\eta_0}$ the coordinate
vectors $\partial / \partial \hat{x}^j$ are Jacobi fields, and are defined even at the conjugate
points.

Finally, we introduce frames $\{ F_j^{t_0}(\hat{x},r) \}_{j=1}^n$
defined by parallel translation along $\gamma^{t_0}_{\hat{x}}$ such that for $j =1,\ ... \ ,\ n$
\[
   F^{t_0}_j(\hat{x},0) = \left.
      \frac{\partial}{\partial \hat{x}^j} \right |_{(\hat{x},0)}
\]
where we are using the notation $r = \hat{x}^n$ and so
\[
   F^{t_0}_n(\hat{x},0) = \dot{\gamma}^{t_0}_{\hat{x}}(0)
\]
points in the opposite direction of $\nu_{0,t_0}(\Psi_{t_0}(\hat{x},0))$. To simplify the presentation we adopt this notation, $r = \hat{x}^n$, throughout the paper. Also we write $\{ f^j_{t_0}(\hat{x},r)\}_{j=1}^n$ for the corresponding dual frame; that is
\[
   \bra f^j_{t_0}(\widehat{x},r),F_k^{t_0}(\widehat{x},r) \cet = \delta_k^j .
\]
where $\langle .,. \rangle$ denotes the usual pairing of $T_{\gamma^{t_0}_{\hat{x}}(r)} \tilde M$ and $T_{\gamma^{t_0}_{\hat{x}}(r)}^* \tilde M$.
In the sequel will also consider the shape operators $S_{r,t}$ when $x_0$ is replaced in the above construction by another point in $\Sigma_{0,t_0}$ represented in the coordinates by $\hat{x}$. We thus have for each $\hat{x}$ and $0 \leq r < t \leq t_0$ such that $r$ and $t$ are not conjugate along $\gamma^{t_0}_{\hat{x}}$ a tensor $S_{r,t}^{t_0}(\hat{x}) \in (T_1^1)_{\gamma^{t_0}_{\hat{x}}(r)} \wt{M}$. We represent $S_{r,t}^{t_0}(\hat{x})$ using the frames constructed above as
\begin{equation}\label{shape_coord}
   S^{t_0}_{r,t}(\hat{x}) = {\bf s}^k_j(\hat{x},t_0;r,t)  f^j_{t_0}(\hat{x},r) \otimes F^{t_0}_k(\hat{x},r).
\end{equation}
For fixed $t_0$ and $\hat{x}$ we will also use the notation $\mathbf{s}^k_j(r,t) = \mathbf{s}^k_j(0,t_0;r,t)$. Note that immediately from the definition we have $\mathbf{s}^{n}_j(\hat{x},t_0;r,t)=\mathbf{s}^{j}_n(\hat{x},t_0;r,t) = 0 $ for all $j$ and because of this in what follows when we write $\mathbf{s}(\hat{x},t_0;r,t)$ (respectively $\mathbf{s}(r,t)$) without indices we will actually be referring to the $(n-1)\times (n-1)$ matrix $\mathbf{s}_{j}^k(\hat{x},t_0;r,t)$ (respectively  $\mathbf{s}_j^k(r,t)$) with $j$, $k = 1,\, ... \, , n-1$. The data for our recovery are the matrix elements ${\bf s}^k_j(\hat{x},t_0;0,t)$, and their first three derivatives with respect to $t$, for $0 < t < t_0$ and $t$ not conjugate to $0$ along $\gamma^{t_0}_{\hat{x}}$.

In~\cite{companion}, we obtain the following result: It is possible to uniquely determine the Riemannian metric $g$ in a neighborhood
of $\gamma_{x_0,\eta_0}([0,t_0))$ in Riemannian normal coordinates
having origin at the point $y_{t_0}$ (this can be done for general metrics, not just ones which are conformally Euclidean). Here, we cast this result
into an algorithm, and construct a conversion from the mentioned
coordinates to Cartesian coordinates, which is the main contribution
of this paper. Essentially, we generalize the time-to-depth conversion
in Dix' original method to multi-dimensional manifolds with Riemannian
metrics conformal to the Euclidean metric and, roughly speaking, show that
if we measure near the point $x_0$ the shape operators of the wave fronts of waves diffracted from the points 
$\gamma_{x_0,\eta_0}(t_0)$, we  can then determine the wave speed $v$
near the geodesic 
$\gamma_{x_0,\eta_0}$. The theoretical contributions of this work are thus contained in the following two theorems.
\medskip\medskip

\begin{theorem}\label{thm1}
Suppose that $M \subset \widetilde{M} = \mathbb{R}^n$ and $\tilde{M} \setminus M$ is open and nonempty, $v \in C^\infty(\widetilde{M})$, $W$ is a set of the form constructed above (see (\ref{Wdef})) using the metric $g = v(x)^{-2} \delta_{ij} \dd x^i \dd x^j$, and $\Sigma_{0,t_0} \subset \widetilde{M} \setminus M$ for all $t_0 \in (0,L)$.  We also assume that $(\widetilde{M},g)$ is complete. If $v |_{(\widetilde{M} \setminus M)\cap W}$ is known then from the data $\mathbf{s}(\hat{x},t_0;0,t)$ (see (\ref{shape_coord})) for all $t_0$ in an open subinterval of $I = (T_1,T_2) \subset (0,L) \cap \mathcal{C}_0(x_0,\eta_0)^c$, $\hat{x} \in U$, and $t \in (0,t_0)\cap \mathcal{C}_{t_0}(\Phi_{t_0}^{-1}(\hat{x}),-\nu_{0,t_0}(\Phi_{t_0}^{-1}(\hat{x})))^c$ it is possible to recover $v$ in a neighborhood of $\gamma_{x_0,\eta_0}([0,T_2))$ in Cartesian coordinates. In dimension strictly larger than $2$ the reconstruction only involves solving ordinary differential equations along the rays of $g$. 
\end{theorem}
\medskip\medskip

Note that the sets $\mathcal{C}_r(x,\eta)$ are always discrete, and so the hypotheses of the theorem assume we have data except for at a discrete set of times. We also stress here that the reconstruction is local along geodesics. That is, to reconstruct $v$ in a neighborhood of $\gamma_{x_0,\eta_0}$ we only require measurements of the shape operator near the point $x_0$ for generalized spheres $\Sigma_{0,t_0}$ centered at points $\gamma_{x_0,\eta_0}(t_0)$ with radii $t_0>0$. We recall $\mathbf{s}(\hat x,t_0;0,t)$ is
representation of the shape operator of these generalized spheres  
in the local coordinates and that we have $x_0\in \Sigma_{0,t_0}$. In the case that the dimension $n \geq 3$ the conversion from Riemannian normal to Cartesian coordinates involves
solving a system of $n + 3 n ^2 + n^3$ nonlinear ordinary differential
equations. In the two dimensional case the construction requires the solution of a Cauchy problem for an elliptic operator (i.e. we must solve the scalar curvature equation \re{scal} with known boundary data). The discretization of the system is directly related to
the available ``density'' of scatterers.

In order to state our second result, we must introduce generalized distance functions. As $\tilde M$ is complete, the map $\exp_y :\ T_y \tilde M \to
\tilde M$ is surjective, and by Sard's theorem, the set ${\cal C}(y)
\subset \tilde M$ of critical values of $\exp_y$ has measure zero. Suppose
$x \in \tilde M \setminus {\cal C}(y)$ and let $\xi \in T_y \tilde M$
be such that $\exp_y(\xi) = x$. Then $\xi$ has a neighborhood $V_{y,\xi}
\subset T_y \tilde M$ such that $\exp_y :\ V_{y,\xi} \to V_{y,\xi}' = \exp_y(V_{y,\xi})$ is a
diffeomorphism; one says that $\exp_y^{-1} :\ V_{y,\xi}' \to V_{y,\xi}$ is a local
inverse of $\exp_y$ corresponding to $\xi$. In $V_{y,\xi}'$, we define the
function, $\rho(\cdot;y,\xi)$, by
\[
  \rho(z;y,\xi) = | \exp_y^{-1}(z) |_g ,
\]
and call this function the generalized distance or travel time
function from the point $y$ associated to direction $\xi$; the wave
fronts observed from a point source at $y$ give us its level sets. We may now state our second theorem.
\medskip\medskip

\begin{theorem}\label{thm2}
Let $M$ and $\wt{M}$ be as in Theorem~\ref{thm1} and suppose that $\partial M$ is smooth and $W' \subset \widetilde{M}$ is an open set such that for all $y \in W' \cap M^{int}$ there is a non-normalized geodesic $\gamma$ for $g$ such that $\gamma([0,1]) \subset W'$, $\gamma(0) = y$, $\gamma(1) \in \partial M$ and $\dot \gamma(1) \notin T \partial M$. We use the notation $\Gamma = W' \cap \partial M$, and also assume that $v |_{W' \cap (\widetilde{M} \setminus M)}$ is known. Further suppose that $\Lambda$ is an indexing set such that there is a bijective map from $\Lambda$ to
\[
\{ (y,\xi)\in T(W' \cap M) \ : \ V'_{y,\xi} \cap \Gamma \not = \emptyset \},
\]
and label $\rho(\cdot;y,\xi) = \rho_\lambda$ according to this map. Then from knowledge of $\rho_\lambda |_{\Gamma}$ for $\lambda \in \Lambda$ we can recover the wave speed $v$ in $W'$.
\end{theorem}
\medskip\medskip

This generalizes an earlier result which says that the boundary distance functions of a compact Riemmanian manifold given on the whole boundary determine the manifold uniquely, see \cite{Ku5} and \cite[Section 3.8]{KKL}. 

The structure of the rest of the paper is as follows. In section~\ref{sec:Prelim} we go over some background from Riemannian geometry that is necessary. Section~\ref{Step1} reviews the first step of the recovery procedure in which we reconstruct the metric $g$ in the coordinates $(\hat{x},r)$ given by the map $\Psi_{t_0}$. Then sections~\ref{Step2} and~\ref{2d} describe the second step of the recovery procedure in which the geodesics and wave speed are reconstructed in Cartesian coordinates respectively in the case of three and higher dimensions and the case of two dimensions. Finally, the proofs of Theorems~\ref{thm1} and~\ref{thm2} are given in section~\ref{sec:proofs}, and section~\ref{sec:conc} contains some concluding remarks. We also provide detail on the conversion from travel time measurements to measurements of the shape operator of generalized spheres for the case when $\partial M$ is flat and $v$ is constant in a neighborhood of $\partial M$ in~\ref{appendA}.

\section{Preliminaries} \label{sec:Prelim}

We summarize the basic differential equations from Riemannian geometry that we will use. We mention some general references to Riemannian geometry
\cite{Eisen, Petersen-book,Petersen-bulletin}. In this work we will use the conventions from \cite{Petersen-book} for the curvature tensor and related quantities in local coordinates.

\subsection{Geodesics}

We evaluate the geodesics by solving
\begin{equation} \label{gamma}
   \frac{\mathrm{d}^2 \gamma^i}{\mathrm{d}\ t^2} +
        \Gamma^i_{kl}(\gamma(t)) \frac{\mathrm{d} \gamma^k}{\mathrm{d} t}
                   \frac{\mathrm{d} \gamma^l}{\mathrm{d} t} = 0 ,
\end{equation}
where
\[
   \Gamma^i_{kl}(x) = \frac 12 g^{pi} \left(\mstrut{0.4cm}\right.
   \frac{\p g_{k p}}{\p x^l} + \frac{\p g_{l p}}{\p x^k}
       - \frac{\p g_{k l}}{\p x^p} \left.\mstrut{0.4cm}\right) 
\]
are the Christoffel symbols. It is also possible to find the solutions of \re{gamma} using the Hamiltonian flow for the Hamiltonian $H(x,p) = \frac 12 p_j p_k g^{jk}(x)$. Although we will not use the Hamiltonian formulation here, we note that it gives the system
\[
\left \{
\begin{array}{l}
\frac{\mathrm{d}}{\mathrm{d} t} \gamma^i = g^{ij}(\gamma(t)) p_j\\
\frac{\mathrm{d}}{\mathrm{d} t} p_i = -\frac{1}{2} p_j p_k \frac{\partial g^{jk}(\gamma(t))}{\partial x^i}
\end{array}
\right .
\]
for the geodesics. From this, we see that, in terms of seismic ray tracing, the geodesics may be identified with generalized image
rays. 

In our case, assuming isotropy (i.e. that the metric is conformally Euclidean), we have
\[
   \Gamma^l_{qm} =
   -\left( \delta_q^l \delta_m^k
      + \delta_m^l \delta_q^k - \delta_{qm} \delta^{kl} \right)
              \frac{\partial f}{\partial x^k} ,\quad
   f = \log(v).
\]
It is more convenient to work with $f$ than $v$ and we will generally do so throughout the remainder of the paper although it is clear that recovery of $f$ is equivalent to recovery of $v$. To make the notation more concise below we introduce the shorthand
\[
\Theta^{lk}_{qm} =
\delta_q^l \delta_m^k + \delta_m^l \delta_q^k - \delta_{qm}
\delta^{kl}.
\]

\subsection{A frame, parallel transport}

As mentioned above, $F^{t_0}_j(\widehat{x},r)$ denotes the parallel
translation of $\partial/\partial \hat{x}^j$ along $\gamma^{t_0}_{\widehat{x}}$ from
$0$ to $r$. Thus, for every $\widehat{x}$ and $r \geq 0$,
$\{F^{t_0}_1(\widehat{x},r),\ldots,F^{t_0}_n(\widehat{x},r)\}$ forms a basis for
$T_{\gamma^{t_0}_{\widehat{x}}(r)} \widetilde{M}$. The invariant formula for parallel translation is
\[
   \nabla_{\dot{\gamma}_{\widehat{x}}(r)} F^{t_0}_k(\widehat{x},r) = 0.
\]
If we introduce matrices which give the frames $F^{t_0}_j(\hat{x},r)$ in Cartesian coordinates
\begin{equation} \label{Fdef}
   F^{t_0}_j(\hat{x},r) = F^k_j(t_0;\hat{x},r)
      \left. \frac{\partial}{\partial x^k}
                     \right |_{\gamma^{t_0}_{\hat{x}}(r)}
\end{equation}
then the invariant formula implies that these satisfy
\[
   \frac{\p}{\p r} F_j^l
   + (\dot{\gamma}^{t_0}_{\widehat{x}})^k(r) \, \Gamma^l_{km}
                    F_j^m = 0
\]
or
\begin{equation} \label{eq:Fpartr}
   \frac{\p}{\p r} F_j^l(t_0;\widehat{x},r)
   + F^k_n(t_0;\widehat{x},r) \, \Gamma^l_{km}(\gamma^{t_0}_{\widehat{x}}(r))
                    F_j^m(t_0;\widehat{x},r) = 0 ,
\end{equation}
and since the coordinates $\hat{x}$ on $\Sigma_{0,t_0}$ are known, we
also know the initial conditions, $F_j^l(t_0;\hat{x},0)$.

\subsection{Curvature}

The Riemannian curvature tensor in any coordinate system is given as a $(1,3)$ tensor field by
\[
   R_{ j kl}^p = \frac{\p }{\p x^j } \Gamma^p _{k l}
        - \frac{\p }{\p x^k} \Gamma^p _{j l}
   + \Gamma^i_{k l} \Gamma^p_{j i} - \Gamma^i_{j l} \Gamma^p_{k i},
\]
or as a $(0,4)$ tensor field as
\[
   \quad g_{ip} R^i_{j kl} = R_{jklp} .
\]
The Ricci curvature tensor is given by the trace
\[
Ric_{i j } = R^k_{k i j }
\]
and the scalar curvature is $scal =
g^{ij} Ric_{i j }$. In the Cartesian coordinates the Riemmanian and Ricci curvature tensors are given respectively by the following formulae in terms of $f$:
\[
R = e^{-2f} \left ( \delta \, \odot\, \left ( - \mathrm{Hess}(f) - \nabla f \cdot \nabla  f  + \frac{1}{2} |\nabla f |^2 \delta \right ) \right )
\]
where $\odot$ is the Kulkarni-Nomizu product (in coordinates, see \cite{Petersen-book} for the invariant formula), $\mathrm{Hess}(f)$ is the Hessian matrix of second derivatives, $\nabla f$ is the (Euclidean) gradient, and $|\nabla f |$ is the Euclidean norm of $\nabla f$; 
\begin{equation} \label{Ric0}
Ric = (n-2) \left ( \mathrm{Hess}(f) +\nabla f \cdot \nabla  f  \right ) + \delta_{ij} \left ( \Delta f + (2-n) \left | \nabla f \right |^2 \right ).
\end{equation}
Also, when the dimension is $n=2$ the scalar curvature satisfies the so-called scalar curvature equation
\begin{equation}\label{scal}
\Delta_g f = \frac{1}{2}\, scal
\end{equation}
where $\Delta_g$ is the Laplace-Beltrami operator corresponding to $g$ given in any coordinate system $\{ y^j \}_{j=1}^2$ by
\[
\Delta_g =  g(y)^{-1/2}
    \frac{\partial}{\partial y^j} g(y)^{1/2}g^{jk}(y)
             \frac{\partial}{\partial y^j}.
\]
Here $g=\hbox{det}([g_{jk}]_{j,k=1}^2 )$ and $[g^{jk}]_{j,k=1}^2$ is
the inverse matrix of $[g_{jk}]_{j,k=1}^2$. From these formulae we see a fundamental difference between the two dimensional case and the case of three or more dimensions. In two dimensions the Ricci curvature and scalar curvature give only the Laplacian of $f$, and so to find $f$ from these curvature tensors would require the solution of an elliptic equation. On the other hand, in three or more dimensions the Ricci curvature tensor depends on all the second partial derivatives of $f$, and in general we can recover a formula for $\mathrm{Hess}(f)$ in terms of the Ricci curvature and $\nabla f$. Indeed, if we define
\begin{equation} \label{Ric}
G = Ric - (n-2) \left (\nabla f \cdot \nabla  f - \left | \nabla f \right |^2 \delta \right )
\end{equation}
then from \re{Ric0} we may calculate
\begin{equation} \label{Ric2}
\mathrm{Hess}(f) = \frac{1}{n-2}\left ( G - \frac{1}{2(n-1)} \mathrm{tr} \left ( G \right ) \delta\right ).
\end{equation}
This is possible in three dimensions, but not in two dimensions, and is the reason we must consider the two cases separately.

We will also write the Riemannian curvature on the geodesic
$\gamma^{t_0}_{\widehat{x}}$ in the frame obtained by parallel transport as
\[
   \wh{R}_{jkl}^p(t_0;\widehat{x},r) F_p^{t_0}(\widehat{x},r)
     = R_{\gamma^{t_0}_{\widehat{x}}(r)}(F_j^{t_0}(\widehat{x},r),
                   F_k^{t_0}(\widehat{x},r)) F_l^{t_0}(\widehat{x},r) ,
\]
or
\[
   \wh{R}_{jkl}^p(t_0;\widehat{x},r)=
   \langle f^p_{t_0}(\widehat{x},r), R_{\gamma^{t_0}_{\widehat{x}}(r)}(F^{t_0}_j(\widehat{x},r),
           F_k^{t_0}(\widehat{x},r))
                F_l^{t_0}(\widehat{x},r) \rangle .
\]
Recalling that $F^{t_0}_n(\widehat{x},r) = \dot \gamma^{t_0}_{\widehat{x}}(r)$, we also
write
\[
   {\bf r}^p_j(t_0;\widehat{x},r) = \wh{R}_{j n n}^p(t_0;\widehat{x},r) ,
\]
and for fixed $t_0$
\begin{equation}
   {\bf r}^p_j(r) := {\bf r}^p_j(t_0;0,r) = \left \langle f_{t_0}^p(0,r) ,
   R_{\gamma_{x_0,\eta_0}(r)}(F^{t_0}_j(0,r), \dot \gamma_{x_0,\eta_0}(r))
              \dot \gamma_{x_0,\eta_0}(r) \right \rangle ,
\end{equation}
for the directional curvature operator which we reconstruct
in the first step of our procedure. Note that as with $\mathbf{s}$, for any $j$ $\mathbf{r}_j^n(r) = \mathbf{r}^j_n(r) = 0$, and so when we write $\mathbf{r}$ without indices we will actually be referring to the corresponding $(n-1) \times (n-1)$ matrix.

We continue in the next sections to describe the actual reconstruction algorithm.

\section{Reconstruction procedure -- Step 1: Determination of the
         metric in $(\hat{x},r)$ coordinates} \label{Step1}

The reconstruction procedure consists of two steps. In the first step
we consider only the single geodesic $\gamma_{x_0,\eta_0}$ and reconstruct $\mathbf{r}^p_j(r)$, and then
the metric $g$ as a function of $r$ for $\hat{x} =0$. In the second step, we
determine $f$, and therefore also $v$, in a neighborhood of $\gamma_{x_0,\eta_0}$ by also varying $\hat{x}$ and $t_0$.

Following the geometric analysis of~\cite{companion} we now describe the first step of the procedure which is itself broken up into a number of substeps below.

\medskip\medskip

\begin{enumerate}
\item
Let $V^j = V^j(r,t)$, $j = 0,\ldots,3$ represent $(n-1) \times (n-1)$
matrices. We solve the autonomous system of ordinary differential
equations for $V(r,t) = (V^j(r,t))_{j=0}^3$,
\begin{equation} \label{EQ A}
\scriptsize
\begin{array}{lcl}
   \pdpd{}{r} V^0 &=& -I - \frac 12 V^0(\T V^3)V^0 ,
   \vspace{.1cm}
\\
\vspace{.1cm}
   \pdpd{}{r} V^1 &=& -\frac 12 (V^1(\T V^3)V^0 + V^0(\T V^3)V^1) ,\\
   \vspace{.1cm}
   \pdpd{}{r} V^2 &=& -\frac 12 (V^2(\T V^3)V^0 + V^0(\T V^3)V^2
                       + 2V^1(\T V^3)V^1) ,\\
                       \vspace{.1cm}
   \pdpd{}{r} V^3 &=& -\frac 12 (V^3(\T V^3)V^0 + V^0(\T V^3)V^3
                       + 3V^2(\T V^3)V^1 + 3V^1(\T V^3)V^2) ,
\end{array}
\normalsize
\end{equation}
in which
\[
   (\T V^3)(r) = V^3(r,r),
\]
for $0 \le r \le t \le t_0$. This system is supplemented with initial
data,
\begin{equation} \label{final init. data}
   V^j(0,t) = \{\p_t^j ({\bf s}(0,t))^{-1}\}_{j=0}^3.
\end{equation}
Applying Picard's theorem in a standard
way, we see that  the equations
(12)-(13) have unique solution on some interval
$t\in [0,t_1]$ (for details on this see~\cite{companion}). In practice, we may use a Runge-Kutta method to solve the system numerically for
$0\leq r \leq t \leq t_1$. The system will not generally have a solution all
the way up to $t_0$; in this case, we must divide the interval
$[0,t_0]$ into several subintervals, and reconstruct on each of these
in turn as described below in substep 3.
\item
We extract the directional curvature operator,
\begin{equation} \label{EQ Ainit}
   {\bf r}(r) = \frac{1}{2} (\T V)(r) .
\end{equation}
Note that this matrix $\mathbf{r}(r)$ is $(n-1) \times (n-1)$, but recall that from this we can recover the full directional curvature operator $\mathbf{r}^j_p(r)$ since the $n$th row and column of $\mathbf{r}^j_p(r)$ are both equal to zero.

\item
\begin{figure} \label{Step1fig}
\centering
\includegraphics{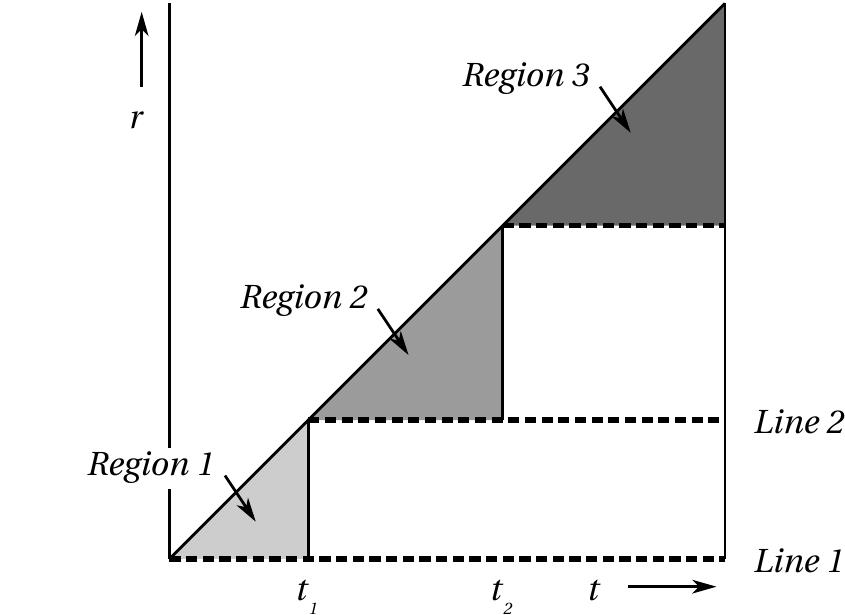}
\caption{A depiction of how the first step of the algorithm proceeds. The original data give $\mathbf{s}$ on line 1. Then substep 1 recovers $\{V^j\}_{j=1}^3$ in region 1 by solving \re{EQ A}. Next substep 2 gives $\mathbf{r}(r)$ for $0 \leq r \leq t_1$ by \re{EQ Ainit}, and then substep 3 gives $\mathbf{s}$ on line 2. Returning to substep 1 again gives $\{V^j\}_{j=1}^3$ in region 2. Then as before substep 2 gives $\mathbf{r}(r)$ for $t_1\leq r \leq t_2$, and substep 3 gives $\mathbf{s}$ on line 3. Continuing in this way we reconstruct $\mathbf{r}$ as far along $\gamma_{x_0,\eta_0}$ as we like.}
\end{figure}
In general the first two steps only reconstruct $\mathbf{r}^j_k(r)$
for $0 \leq r \leq t_1$ where $t_1 < t_0$ since we may not able to
solve (\ref{EQ A}) all the way up to $t_0$. Here, we describe how to
find $\mathbf{r}^j_k(r)$ for $r$ in the entire interval $[0,t_0]$. The
idea is to replace $0$ by $t_1$, and then use knowledge of
$\mathbf{s}_j^k(t_1,t)$ for $t> t_1$. Indeed, let us fix $t >
t_1$. Now we find the coordinates of the Jacobi fields corresponding to the coordinate vectors fields of $(\hat{x},r)$ along $\gamma_{x_0,\eta_0}$ with respect to the parallel
frame. These can be found by solving the system
\begin{equation}
   \pdpd{}{r} \left(\begin{array}{r}
        {\bf j}^j_{k}(r,t) \\
        \dot {\bf j}^j_{k}(r,t) \end{array}\right)
   = \left(\begin{array}{cc} 0 & \delta^j_p \\
               -{\bf r}^j_p(r) & 0 \end{array}\right)
     \left(\begin{array}{r}  {\bf j}^p_{k}(r,t) \\
        \dot {\bf j}^p_{k}(r,t) \end{array}\right) ,
\end{equation}
with initial conditions
\begin{equation}
   \left(\begin{array}{r}
        {\bf j}^j_{k}(0,t) \\
        \dot {\bf j}^j_{k}(0,t) \end{array}\right)
   = \left(\begin{array}{c}
        \delta^j_k \\ -\mathbf{s}_k^j(0,t) \end{array}\right) .
\end{equation}
We can then recover $\mathbf{s}(t_1,t)$ by the equation
\[
   \pdpd{\mathbf{j}}{r}(t_1,t) \,
       (\mathbf{j}^{-1})(t_1,t) = -\mathbf{s}(t_1,t) .
\]
Now we may return to substep 1 and solve (\ref{EQ A}) with (\ref{final init. data})
replaced by
\[
   V^j(t_1,t) = \{\p_t^j ({\bf s}(t_1,t))^{-1}\}_{j=0}^3 .
\]
This then allows us to use \re{EQ Ainit} to recover $\mathbf{r}(r)$ for $r$
up another time $t_2$ say with $t_1<t_2\leq t_0$. If $t_2$ is less
than $t_0$ we repeat this same procedure again with $t_1$ replaced by
$t_2$, and so on. By results in~\cite{companion} there is a lower bound on the size of each step we take (i.e. a lower bound on $t_i - t_{i-1}$), and so by induction we eventually recover
$\mathbf{r}^j_k(r)$ for $r$ on the entire interval $[0,t_0]$. For a visual depiction of how the first three substeps proceed see Figure~\ref{Step1fig}.

\item
We now obtain the metric in the $(\hat{x},r)$ coordinates given by the coordinate map $\Psi_{t_0}$ along $\gamma_{x_0,\eta_0}$, which we write as $\wh{g}_{jk}(0,r)$, by the formula
\begin{eqnarray}\label{hatg}
   \wh{g}_{jk}(0,r)
   = {\bf j}^p_{j}(r,t_0)
       {\bf j}^q_{k}(r,t_0) \, \mathring{g}_{pq} ,
\end{eqnarray}
where
\[
\mathring{g}_{pq} = g(F_p(0,0),F_q(0,0)) = v(x_0) F_p^j(0,0) F_q^i(0,0) \delta_{ij}.
\]
We assume in the hypotheses of Theorem~\ref{thm1} that $v(x_0)$ is known, and $F_p^j(0,0)$ (resp. $F_q^j(0,0)$) are the components of the coordinate vector $\partial/\partial \hat{x}^p$ (resp. $\partial/ \partial \hat{x}^q$) with respect to Cartesian coordinates. Thus $\mathring{g}_{pq}$ is known. By adjusting the choice of $x_0$ we also find the metric with respect to the $(\hat{x},r)$ coordinates where they are defined (i.e. on all of $W(t_0)$).

\end{enumerate}

\section{Reconstruction procedure -- Step 2: Transformation of
         coordinates} \label{Step2}

By the procedure described in the previous section we can reconstruct the metric $\hat{g}_{jk}(\hat{x},r)$ in the coordinates $(\hat{x},r)$ everywhere in the domain $W(t_0)$ where those coordinates are defined. Thus we can also reconstruct the Ricci curvature tensor in these coordinates and also the scalar curvature, which we will use below. In this section, we show how to determine the velocity function $v(\gamma^{t_0}_{\widehat{x}}(r))$ and the geodesics $\gamma^{t_0}_{\widehat{x}}(r)$ in the Cartesian coordinates from this information. The reconstruction is done initially up to the first conjugate point to $t_0$ along $\gamma_{x_0,\eta_0}$. Then, as described at the end of this section, $t_0$ must be changed to allow reconstruction past this conjugate point.

As observed above, the coordinate vectors in the $(\hat{x},r)$ coordinates are Jacobi fields along $\gamma^{t_0}_{\hat{x}}$; in particular, if
\[
   \left. \frac{\partial}{\partial \hat{x}^k}
       \right|_{\gamma^{t_0}_{\hat{x}}(r)}
              = \mathbf{j}^l_k(t_0;\hat{x},r) F^{t_0}_l(\hat{x},r) = \mathbf{j}^l_k(t_0;\hat{x},r) F_l^p(t_0;\hat{x},r) \left .\frac{\partial}{\partial x^p} \right |_{\gamma^{t_0}_{\hat{x}}(r)},
\]
(recall that $\frac{\partial}{\partial x^p}$ are the coordinate vectors for the Euclidean coordinates) then the matrix $\mathbf{j}^l_k(t_0;\hat{x},r)$ satisfies
\begin{equation}\label{Jacobisys}
   \pdpd{}{r} \left(\begin{array}{r}
        {\bf j}^l_{k}(t_0;\widehat{x},r) \\
        \dot {\bf j}^l_{k}(t_0;\widehat{x},r) \end{array}\right)
   = \left(\begin{array}{cc} 0 & \delta^l_p \\
            -{\bf r}^l_p(t_0;\widehat{x},r) & 0 \end{array}\right)
     \left(\begin{array}{r} {\bf j}^p_{k}(t_0;\widehat{x},r) \\
        \dot {\bf j}^p_{k}(t_0;\widehat{x},r) \end{array}\right) ,
\end{equation}
supplemented with the initial data
\[
   \left.\left(\begin{array}{r}
        {\bf j}^l_{k}(t_0;\widehat{x},0) \\
        \dot {\bf j}^l_{k}(t_0;\widehat{x},0) \end{array}\right)
   \right|_{r=0} = \left(\begin{array}{c}
        \delta^l_k \\ -\mathbf{s}^l_k(\hat{x},t_0;0,t_0) \end{array}\right) .
\]
In fact, here we are simply adding the dependence on $\hat{x}$ and $t_0$ to the
same quantities already considered in the previous section. Since parallel translation preserves the
metric, we also have the relation
\[
   v(\gamma^{t_0}_{\hat{x}}(0))^{-2} F^l_j(t_0;\widehat{x},0)
            F^q_k(t_0;\widehat{x},0) \delta_{lq}
   = v(\gamma^{t_0}_{\widehat{x}}(r))^{-2}
        F^l_j(t_0;\widehat{x},r) F^q_k(t_0;\widehat{x},r) \delta_{lq}.
\]
By taking determinants of the matrices on each side and then the natural log, we obtain the following formula
\begin{equation}\label{vrel}
   f(\gamma^{t_0}_{\widehat{x}}(r)) =\frac{1}{n} \log \left(
      v(\gamma^{t_0}_{\hat{x}}(0))^n \left|\frac{ \mathrm{det}(F^l_j(t_0;\widehat{x},r))}{\mathrm{det}(F^l_j(t_0;\widehat{x},0))}
          \right| \right).
\end{equation}
Recall that $F^l_j(t_0;\hat{x},r)$ are the components of the parallel frame with respect to the Cartesian coordinate vectors (cf. (\ref{Fdef})). Also, combining some of the previous formulas we have
\begin{equation}\label{coc1}
   \frac{\partial}{\partial \widehat{x}^k} f(\gamma^{t_0}_{\widehat{x}}(r))
   = {\bf j}_{k}^p(t_0;\widehat{x},r) F_p^l(t_0;\widehat{x},r)
            \frac{\partial f}{\partial x^l} (\gamma^{t_0}_{\widehat{x}}(r)).
\end{equation}
Away from conjugate points along $\gamma^{t_0}_{\widehat{x}}$, we may invert
to obtain
\begin{equation}\label{coc}
   \frac{\partial f}{\partial x^l} (\gamma^{t_0}_{\widehat{x}}(r))
   = (F^{-1})^p_l(t_0;\wh{x},r) ({\bf j}^{-1})^j_p(t_0;\wh{x},r)
    \frac{\partial}{\partial \widehat{x}^j}
                f(\gamma^{t_0}_{\widehat{x}}(r)) 
\end{equation}
which may further be combined with \re{vrel} to show
\small
\begin{equation} \label{dfdx}
\frac{\partial f}{\partial x^l} (\gamma^{t_0}_{\widehat{x}}(r)) =\frac{1}{n} (F^{-1})^p_l(t_0;\wh{x},r) ({\bf j}^{-1})^j_p(t_0;\wh{x},r) \left (\frac{\partial \sigma^{t_0}}{\partial \hat x^j}(\hat{x}) + (F^{-1})_b^c(t_0;\hat{x},r) \frac{\partial F_c^b}{\partial \hat{x}^j}(t_0;\hat{x},r) \right )
\end{equation}
\normalsize
where we use the notation
\[
   \sigma^{t_0}(\wh{x}) := \log\left( \frac{v(\gamma^{t_0}_{\hat{x}}(0))^n}{
   \left| \mathrm{det}(F^l_j(t_0;\widehat{x},0)) \right|} \right) .
\]
\begin{remark} \label{rem}
{\rm A ``stepping" recovery procedure in the spirit of Dix' original method may now be presented. We introduce
a step size $h$ in $r$, and for $\alpha \in \mathbb{N}$ we label
$r_\alpha = \alpha h$. If we know
$f(\gamma^{t_0}_{\widehat{x}}(r_{\alpha-1}))$,
$F_j^l(t_0;\widehat{x},r_{\alpha-1})$, and ${\bf
  j}_k^l(t_0;\widehat{x},r_{\alpha-1})$, then we approximate the same
quantities at $r_\alpha$ by the following strategy. First we
numerically estimate the derivatives
\[
   \frac{\partial}{\partial \wh{x}^k}
          f(\gamma^{t_0}_{\wh{x}}(r_{\alpha-1})) .
\]
Then we use (\ref{coc}) to estimate the derivatives
\[
   \frac{\partial f}{\partial x^k} (\gamma^{t_0}_{\wh{x}}(r_{\alpha-1})) .
\]
We then have an estimate of
$\Gamma^l_{km}(\gamma^{t_0}_{\widehat{x}}(r_{\alpha - 1}))$. Next, we use
this estimate to perform a forward Euler step in (\ref{eq:Fpartr}) and
get an approximation of $F^l_j(t_0;\wh{x},r_\alpha)$. Then, finally, we
use (\ref{vrel}) to obtain the approximation for
$f(\gamma^{t_0}_{\wh{x}}(r_\alpha))$. We note that ${\bf
  j}_k^l(t_0;\widehat{x},r_\alpha)$ may be obtained through a completely
independent calculation using (\ref{Jacobisys}). Because
$F^l_n(t_0;\wh{x},r) = (\dot \gamma_{\wh{x}}^{t_0})^l(r)$, we can then approximate
$(\gamma_{\wh{x}}^{t_0})^l(r_\alpha)$ as well.}
\end{remark}

\subsection*{A closed system of ordinary differential equations for
             $n \ge 3$}

While the technique described in remark~\ref{rem} may be a direct generalization of Dix' method, in fact we seek a
single closed system of ordinary differential equations that could be solved
using any numerical scheme to give all the desired quantities. This is
possible, although the method we present here works in three or more
dimensions only. The reason for this limitation, as explained above, comes from our use of formulas \re{Ric} and \re{Ric2} to express the Hessian of $f$ in terms of the Ricci curvature and the first order derivatives of $f$ in Cartesian coordinates. Since we actually know the Ricci curvature in $(\hat{x},r)$ coordinates given by $\Psi_{t_0}$, which we will label as $\wh{Ric}_{pq}$, we also need the formula for the tensorial change
\begin{eqnarray} \label{CartRic}
Ric_{ij}(\gamma^{t_0}_{\hat{x}}(r)) = && \\
&&\hskip-1in \wh{Ric}_{pq}(t_0;\hat{x},r) (\mathbf{j}^{-1})_l^p(t_0;\hat{x},r) (F^{-1})_i^l(t_0;\hat{x},r) (\mathbf{j}^{-1})_m^q(t_0;\hat{x},r) (F^{-1})_j^m(t_0;\hat{x},r). \nonumber
\end{eqnarray}
Now we can describe how to get the closed system of ordinary differential equations.

We differentiate (\ref{eq:Fpartr}) and use (\ref{coc1}) to get
\begin{equation}\label{dFsys}
   \frac{\partial}{\partial r}
   \frac{\partial F^l_j}{\partial \wh{x}^a}
   =  \frac{\partial F^q_n}{\partial \wh{x}^a}
          F_j^m \Theta^{lk}_{qm} \frac{\partial f}{\partial x^k}
     + F^q_n \frac{\partial F_j^m}{\partial \wh{x}^a}
         \Theta^{lk}_{qm} \frac{\partial f}{\partial x^k}
     + F^q_n F_j^m \Theta^{lk}_{qm} 
         {\bf j}_a^p F_p^c
                \frac{\partial^2 f}{\partial x^k \partial x^c} .
\end{equation}
where we have  suppressed the dependence on
$(t_0;\wh{x},r)$. We may use (\ref{Ric}), (\ref{Ric2}), \re{dfdx}, and \re{CartRic} to express the right-hand side
of (\ref{dFsys}) only in terms of
$\mathbf{j}$, $F$, derivatives of $F$, and the Ricci curvature $\widehat{Ric}_{pq}$. Combining all the previous equations we now have a closed system of ordinary differential equations that may be solved uniquely up to conjugate points. In the next paragraph, we
summarize the entire method for convenience.

As claimed above we have now produced a closed system of ordinary
differential equations that may be solved to obtain $v$ and
$\gamma^{t_0}_{\wh{x}}$ in Cartesian coordinates. The system is nonlinear and contains $n+ 3n^2 +
n^3$ equations. It may be written as
\begin{equation}\label{W sys}
   \frac{\partial}{\partial r}
   \left(\begin{matrix}
         \gamma^l \\[0.1cm]
         {\bf j}_{k}^l \\[0.1cm]
         \dot {\bf j}_{k}^l \\[0.1cm]
         F^l_k \\[0.1cm]
         \frac{\partial F^l_k}{\partial \wh{x}^p}
         \end{matrix}\right)
   = \left(\begin{matrix}
     W_\gamma^l(r,\wh{x},{\bf j},\dot {\bf j}, 
                    F,\frac{\partial F}{\partial \wh{x}}) \\[0.1cm]
     W_{{\bf j};k}^l(r,\wh{x},{\bf j},\dot {\bf j}, 
                    F,\frac{\partial F}{\partial \wh{x}}) \\[0.1cm]
     W_{\dot {\bf j};k}^l(r,\wh{x},{\bf j},\dot {\bf j}, 
                    F,\frac{\partial F}{\partial \wh{x}}) \\[0.1cm]
     W_{F;k}^l(r,\wh{x},{\bf j},\dot {\bf j}, 
                    F,\frac{\partial F}{\partial \wh{x}}) \\[0.1cm]
     W_{\frac{\partial F}{\partial \wh{x}};kp}^l(r,\wh{x},{\bf j},
     \dot {\bf j},F,\frac{\partial F}{\partial \wh{x}})
     \end{matrix}\right) .
\end{equation}
We describe how each of the ``$W$'' functions on the right-hand side are to be evaluated. $W_\gamma^l$ and $W_{{\bf j};k}^l$ are the
simplest. Recalling that $F^{t_0}_n(\hat{x},r) = \dot \gamma^{t_0}_{\hat{x}}(r)$ and using \re{Jacobisys} they are given by
\[
   W_\gamma^l = F_n^l \quad \mbox{and} \quad
   W_{{\bf j};k}^l = \dot {\bf j}_{k}^l .
\]
Next, again according to (\ref{Jacobisys}), $W_{\dot {\bf j};k}^l$ is given
by
\[
   W_{\dot {\bf j};k}^l = -{\bf r}^l_j {\bf j}_{k}^j .
\]
Since $\wh{x}^n = r$, $W_{F;k}^l$ is given by
\[
   W_{F;k}^l = \frac{\partial F^l_k}{\partial \wh{x}^n} .
\]
Finally, $W_{\frac{\partial F}{\partial \wh{x}};kp}^l$ is given by
(\ref{dFsys}) where we calculate $\pdpd{f}{x^j}$ using (\ref{dfdx}), and calculate
$\frac{\partial^2 f}{\partial x^i \partial x^j}$ in several steps using the values of $\pdpd{f}{x^j}$ already calculated, \re{CartRic}, (\ref{Ric}), and then (\ref{Ric2}). System (\ref{W sys}) gives a nonlinear system of ODEs which can be solved to recover $\gamma$ and $F_k^l$ up to conjugate points. However, at the conjugate points the matrix $\mathbf{j}$ will become singular, and we will not be able to continue.

Note that we do not explicitly solve for $v$ in this system, but
after we have found $F^l_k$, then $f$ and therefore $v$ may be
calculated from (\ref{vrel}). Indeed since $f = \mathrm{log}(v)$, (\ref{vrel}) becomes
\[
v(\gamma^{t_0}_{\hat{x}}(r)) = v(\gamma^{t_0}_{\hat{x}}(0)) \left | \frac{\mathrm{det}(F^l_j(t_0;\hat{x},r))}{\mathrm{det}(F^l_j(t_0;\hat{x},0))} \right |^{1/n}.
\]

\medskip\medskip

\subsection*{The presence of conjugate points} \label{Sec conjugate}

\begin{figure}
\center{
\includegraphics[width=.75\textwidth]{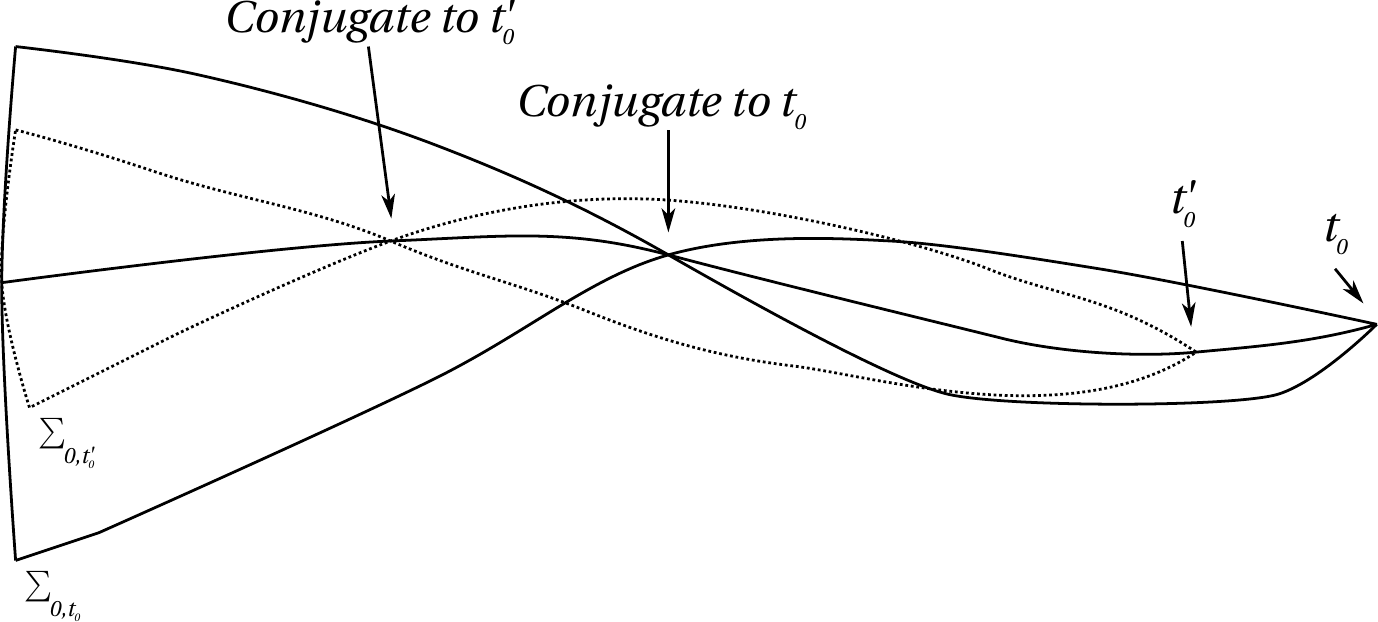}}
\caption{When there are conjugate points we may have to perform step 1 again with $t_0$ replaced by another value $t_0'$.}
\end{figure}

By its construction, the system (\ref{W sys}) has a solution up to the first conjugate point of $t_0$ along
$\gamma_{x_0,\eta_0}$, and since the functions on the right hand side of (\ref{W sys}) are Lipschitz continuous away from the conjugate points this solution is unique. Thus we can recover $v$ along the entire length of $\gamma_{x_0,\eta_0}$ using (\ref{W sys}) if there is no
conjugate point to $t_0$. However, if there is a conjugate point, then we must
follow a more sophisticated strategy. 

Note that once we are able to calculate the matrix $\mathbf{j}(t_0;\hat{x},r)$ after step 1, then we can identify all of the conjugate points of $t_0$. Since we are only concerned with a finite length along the geodesic $\gamma_{x_0,\eta_0}$, by the Morse Index Theorem (see e.g. \cite[Theorem 15.1]{Milnor}) there are only a finite number of points conjugate to $t_0$, and each of these conjugate points also has only a finite number of conjugate points on the finite interval of interest. Thus we can pick an alternate value of $t_0'$ such that $t_0$ and $t_0'$ do not have any of the same conjugate points, and there are no conjugate points to either on $[t_0',t_0]$. Performing step 1 with $t_0$ replaced by by $t'_0$, we can also consider the system (\ref{W sys}) with $t_0$ replaced by $t'_0$. 

Now suppose that $\{t_{j}\}_{j=1}^m$ is the set of times conjugate to $t_0$ in the interval $[0,t_0)$ given in decreasing order (i.e. so that $t_j > t_{j+1}$ for all $j$). Then for each $j$ there is an open interval $I'_j \subset [0,t_0')$ containing $t_j$ such that $I_j$ does not contain any times conjugate to $t_0'$ along $\gamma_{x_0,\eta_0}$. Suppose that $\{I_j\}_{j=1}^{m+1}$ is another collection of open intervals so that $\{I_j \}_{j=1}^{m+1} \cup \{I_j'\}_{j=1}^m$ covers $[0,t_0)$ and $I_j \subset [ t_{j},t_{j-1} ]$ for all $j$ where we are setting $t_{m+1} = -\epsilon$ for some $\epsilon>0$ sufficiently small. 

Now, let us set $U_{m+1} = U$ (recall that $U$ is the domain of the inverse coordinate map $\Phi_{t_0}$). By solving (\ref{W sys}) and possibly shrinking $U_{m+1}$ we can find the wave speed $v$ in Cartesian coordinates on the set $W_{m+1} := \Psi_{t_0}(U_{m+1} \times I_{m+1})$. By the continuity of $\Psi_{t_0'}$ we can find an open set $U_m' \subset \mathbb{R}^{n-1}$ and an open interval $J'_{m+1}$ containing $0$ and intersecting $I'_m$ nontrivially such $\Psi_{t_0'}(U_m' \times J'_{m+1} ) \subset W_{m+1}$. Therefore, since the wave speed is known in $W_{m+1}$ we can find the parallel fields and Jacobi fields in Cartesian coordinates corresponding to $t_0'$ for $\hat{x} \in U_m'$ and $r \in J_{m+1}'$. Now, after possibly shrinking $U_m'$, we can solve (\ref{W sys}) with $t_0$ replaced by $t_0'$ to find the wave speed in Cartesian coordinates on $\Psi_{t_0'}(U_m' \times I_m')$. Thus we have reconstructed the wave speed on $\Psi_{t_0'}(U_m' \times (I_{m+1} \cup I_m'))$. Switching the roles of $t_0$ and $t_0'$ we can reconstruct the wave speed on a set of the form $\Psi_{t_0}(U_m \times (I_{m+1} \cup I_m' \cup I_m))$ for some open set $U_m \subset U_{m+1}$. Continuing in this way after a finite number of steps we can eventually reconstruct the wave speed on a set of the form $\Psi_{t_0}(U_{1} \times (-\epsilon,t_0))$ which is a neighborhood of the geodesic $\gamma_{x_0,\eta_0}$. This completes the reconstruction in dimension three or higher.

\section{Two-dimensional case}\label{2d}

The method of the previous section, step 2 in the reconstruction procedure, will not work in two dimensions. The basic reason for this is that we cannot determine all of the second partial derivatives of $f$ from the curvature of $g$, but rather can only obtain the Laplacian of $f$ as shown in \re{scal}. Therefore in the two-dimensional case we are forced to recover $f$ by solving \re{scal}. We discuss this in more detail below, but first we will also revisit step 1 in the two-dimensional case where the formulae can be simplified.

The main simplification comes from the fact that in the two dimensional case the trace of $S^{t_0}_{r,t}(\hat{x})$ contains all of the same information as $S^{t_0}_{r,t}(\hat{x})$ itself, and so it is actually easier to consider this trace. Indeed, let us define
\[
\alpha(\hat{x},t_0;r,t) = \mathrm{tr}\left ( S^{t_0}_{r,t}(\hat{x}) \right ).
\]
In this case we have the following simple method of calculating $\alpha(\hat{x},t_0;r,s)$ away from conjugate points. If $r$ is not conjugate to $t$ along $\gamma^{t_0}_{\hat{x}}$, then there is a distance function defined for $(\hat{z},s)$ in a neighborhood of $(\hat{x},r)$ by
\[
d^{t_0}_{\hat{x},r,t}(\hat{z},s) = \left | \mathrm{exp}^{-1}_{\gamma^{t_0}_{\hat{x}}(t)}(\gamma^{t_0}_{\hat{z}}(s)) \right |_g.
\]
In the seismic context this is nothing other than the local travel time along rays close to $\gamma^{t_0}_{\hat{x}}$ from $\gamma^{t_0}_{\hat{x}}(t)$ to points near $\gamma^{t_0}_{\hat{x}}(r)$. By \cite[p.46]{Petersen-book} we have
\[
\alpha(\hat{x},t_0;r,t) = \Delta_g d^{t_0}_{\hat{x},r,t}.
\]
We continue to review step 1 in the two dimensional case.

\subsection{Step 1 redux: The two dimensional case}
We note that the $\{ V^j \}_{j=0}^3$ which appear in equation \re{EQ A} are actually all scalars and the only nonzero component of $\mathbf{r}_j^k$ is $\mathbf{r}_1^1$ which is what we recover in \re{EQ Ainit}. Note also that in the two dimensional case, $\mathbf{r}_1^1$ contains the same information as the sectional curvature (if $F_1$ is chosen to have unit length with respect to $g$ then in fact $\mathbf{r}_1^1$ is the sectional curvature).

The other simplification occurs in substep 4. In the two dimensional case the metric must have the form
\[
   \wh{g}_{jk} = \varphi^{t_0}(\hat{x},r)^2 d\hat{x}^2 + dr^2 ,
\]
Now by \cite[p.46]{Petersen-book}, $\varphi^{t_0}(\hat{x},r)$ satisfies the
equation
\begin{eqnarray} \label{EQ: Jacobi e.}
   \frac{\partial^2}{\partial r^2}
        \varphi^{t_0}(\hat{x},r) + \mathbf{r}_1^1(t_0,\hat{x},r) \varphi^{t_0}(\hat{x},r) = 0 
\end{eqnarray} 
where $\mathbf{r}_1^1(t_0,\hat{x},r)$ is already known. Since the metric is also known in a neighborhood of $\Sigma_{0,t_0}$ we may simply solve this equation with the known initial data $\varphi^{t_0}(\hat{x},0)$ and $\partial_r \varphi^{t_0}(\hat{x},0)$ in order to find the metric in the $(\hat{x},r)$ coordinates defined by $\Psi_{t_0}$. We see that it is not necessary in this case to compute the matrix $\mathbf{j}$ corresponding to the Jacobi fields.

Now we continue to show how step 2 may be accomplished in the two dimensional case. The method makes use of the scalar curvature rather than the Ricci curvature.

\subsection{Step 2 redux: The two dimensional case via the scalar curvature equation}

In step 2 we must take a different strategy for the two dimensional case. The difference is that we cannot express the second derivatives of $f$ which appear in \re{dFsys} in terms of $\mathbf{j}$, $F$, derivatives of $F$, and the Ricci curvature. Instead we use a method inspired by the treatment from \cite[Section 4.5.6]{KKL} of a different problem. After we recover the metric $g$ in the coordinates $(\hat{x},r)$ given by $\Psi_{t_0}$ we use the scalar curvature equation \re{scal} to directly solve for $f$ in these coordinates. Indeed, we assume that 
\begin{equation} \label{data 1a}
   v|_{\Sigma_{0,t_0}}\hbox{ and $\frac{\partial v}{\partial r}|_{\Sigma_{0,t_0}}$}
\end{equation}
are known, and so in fact we have Cauchy data for $f$ on $\Sigma_{0,t_0}$. Thus $f$ satisfies a Cauchy problem for the elliptic operator $\Delta_g$ (see \re{scal}) expressed in $(\hat{x},r)$ coordinates. This equation is given explicitly by
\[
\hat{g}(t_0;\hat{x},r)^{-1/2} \frac{\partial }{\partial \hat{x}^j} \hat{g}(t_0;\hat{x},r)^{1/2} \hat{g}^{jk}(t_0;\hat{x},r) \frac{\partial}{\partial \hat{x}^k} f(\gamma^{t_0}_{\hat{x}}(r)) = \frac{1}{2} \ scal(\gamma^{t_0}_{\hat{x}}(r)).
\]
Here $\hat{g}^{jk}(t_0;\hat{x},r)$ is the inverse of the matrix given by formula (\ref{hatg}) extended to values of $\hat{x}$ other than $0$ and $\hat{g}(t_0;\hat{x},r)$ is the determinant of the matrix $\hat{g}_{jk}(t_0;\hat{x},r)$. The scalar curvature, $scal$, can be computed from $\hat{g}_{jk}(\hat{x},r)$. Expressed in the coordinates this is an elliptic equation although it is degenerate when $r$ is conjugate to $t_0$ along $\gamma^{t_0}_{\hat{x}}$. Adopting the notation of the previous section, by the unique continuation principal (for a modern review of Cauchy problems for elliptic operators see \cite{Aless}) we can first reconstruct $f(\gamma_{\hat{x}}^{t_0}(r))$ for $(\hat{x},r) \in U_{m+1} \times I_{m+1}$. We note, however, that this reconstruction is generally unstable (once again see \cite{Aless} for a detailed review of the stability of this type of problem).

Now once we have recovered $f$ in $(\hat{x},r)$ coordinates for $(\hat{x},r) \in U_{m+1} \times I_{m+1}$ the system \re{W sys} can be replaced by a significantly simpler system. Indeed, if we combine the equation \re{Jacobisys} for the Jacobi field matrix with \re{eq:Fpartr} and \re{coc}, then we have a closed system of ordinary differential equations which may be solved just like \re{W sys} for the higher dimensional case. For convenience we write this system down explicitly. The system is
\begin{equation}\label{W sys2}
   \frac{\partial}{\partial r}
   \left(\begin{matrix}
         \gamma^l \\[0.1cm]
         {\bf j}_{k}^l \\[0.1cm]
         \dot {\bf j}_{k}^l \\[0.1cm]
         F^l_k \\[0.1cm]
         \end{matrix}\right)
   = \left(\begin{matrix}
     W_\gamma^l(r,\wh{x},{\bf j},\dot {\bf j}, 
                    F) \\[0.1cm]
     W_{{\bf j};k}^l(r,\wh{x},{\bf j},\dot {\bf j}, 
                    F) \\[0.1cm]
     W_{\dot {\bf j};k}^l(r,\wh{x},{\bf j},\dot {\bf j}, 
                    F) \\[0.1cm]
     W_{F;k}^l(r,\wh{x},{\bf j},\dot {\bf j}, 
                    F) \\[0.1cm]
     \end{matrix}\right) .
\end{equation}
Here $W_\gamma^l(r,\wh{x},{\bf j},\dot {\bf j}, F)$, $W_{{\bf j};k}^l(r,\wh{x},{\bf j},\dot {\bf j}, F)$, and $W_{\dot {\bf j};k}^l(r,\wh{x},{\bf j},\dot {\bf j},F)$ are given by the same formulae shown below \re{W sys}. The last entry on the right hand side is given, according to \re{eq:Fpartr} and \re{coc}, by
\[
W_{F;k}^l(r,\wh{x},{\bf j},\dot {\bf j}, F) = \Theta^{lj}_{pq} F^p_n F^q_k (F^{-1})^i_j (\mathbf{j}^{-1})_i^a \frac{\partial}{\partial \hat{x}^a} f(\gamma^{t_0}_{\hat{x}}(r)).
\]
Solving these equations we can recover $f$, and therefore also $v$, in Cartesian coordinates on $\Psi_{t_0}(U_{m+1}\times I_{m+1})$. Once we have this we can find Cauchy data on the surface $U'_m \times \{t'_m\}$ for some $t'_m \in I_m' \cap \{ t < t_m \}$ for the scalar curvature equation in the coordinates given by $\Psi_{t_0'}$, and repeat the process described above. Thus we can introduce a stepping procedure as in the higher dimensional case and this completes the reconstruction in the case of two dimensions.

\section{Proofs}\label{sec:proofs}

\proof{\ref{thm1}}
This theorem follows from the recovery procedure that we have presented in sections~\ref{Step1} through~\ref{2d}. Indeed, from results in~\cite{companion} applying step 1, which is described in section~\ref{Step1}, for any $\hat{x} \in U$ we can recover the metric $g$ in $(\hat{x},r)$ coordinates corresponding to $\Psi_{t_0}$ for any $t_0\in I$. Then, in dimension $3$ or higher, the argument of section~\ref{Step2} completes the proof, while in dimension $2$ we must use the scalar curvature equation as described in section~\ref{2d}.
\endproof

\noindent The proof of Theorem~\ref{thm2} involves reduction to the case of Theorem~\ref{thm1}. We also perform this reduction more explicitly in the case that $\partial M$ is flat and $v$ is constant in a neighborhood of $\partial M$ in the appendix.

\proof{\ref{thm2}}

Let us begin by taking any $\lambda_0 \in \Lambda$. By the hypotheses there exists a point $\wt{z} \in \mathrm{domain}(\rho_{\lambda_0}) \cap \Gamma$. Let $\hat{x} = (\hat{x}^1, \ ... \ , \hat{x}^{n-1})$ be a set of local coordinates defined on an open subset of $\mathrm{domain}(\rho_{\lambda_0}) \cap \Gamma$ containing $\wt{z}$ and denote by $\Phi_z$ the coordinate mapping which we suppose has range given by $U \subset \mathbb{R}^{n-1}$. We construct $\Phi_z$ such that $\Phi_z(\wt{z}) = 0$. Finally, let $\nu$ denote the outward pointing unit normal, with respect to $g$, vector field for $\partial M$.

We will now write $\mathrm{d}_{\partial M}$ for the exterior derivative on $\partial M$. For $\hat{x} \in U$ and $\epsilon >0$ let us consider the set
\begin{eqnarray}
\wt{\gamma}_{\hat{x},\epsilon} & = & \bigg \{ \lambda \in \Lambda \ : \ \mathrm{d}_{\partial M}\rho_\lambda ( \Phi^{-1}_z(\hat{x})) = \mathrm{d}_{\partial M}\rho_{\lambda_0} ( \Phi^{-1}_z(\hat{x})), \label{eq:wtgamma}\\
&&\hskip.5in \mbox{and} \ \rho_\lambda(\Phi^{-1}_z(\hat{x})) < \rho_{\lambda_0} ( \Phi^{-1}_z(\hat{x})) +\epsilon \bigg \} \nonumber
\end{eqnarray}
and mapping $\alpha_{\hat{x},\epsilon}: \wt{\gamma}_{\hat{x},\epsilon} \rightarrow (-\epsilon, \rho_{\lambda_0}(\Phi^{-1}_z(\hat{x})))$ defined by
\[
\alpha_{\hat{x},\epsilon}(\lambda) = \rho_{\lambda_0}(\Phi^{-1}_z(\hat{x}))- \rho_\lambda(\Phi^{-1}_z(\hat{x})).
\]
This set and mapping can be found from the data. Now, if $\mathrm{range}(\alpha_{\hat{x},\epsilon})$ is dense in $(-\epsilon,\rho_{\lambda_0}(\Phi_{z}^{-1}(\hat{x}))$ for all $\hat{x}$ in a neighborhood of $0$ and $|\mathrm{d}_{\partial M}\rho_{\lambda_0} ( \Phi^{-1}_z(\hat{x}))|_g < 1$ then the geodesic represented by $\wt{\gamma}_{0,\epsilon}$ (see the next paragraph) is contained in $W'$ and intersects $\partial M$ transversally. In this case we continue. Otherwise we do not use this choice of $\lambda_0$.

Intuitively, $\wt{\gamma}_{\hat{x},\epsilon}$ is a representation of the portion of a geodesic passing through $\Phi_z^{-1}(\hat{x})$ which lies inside $M$, and the range of $\alpha_{\hat{x},\epsilon}$ parametrizes this geodesic segment. Indeed, for every $\lambda \in \wt{\gamma}_{\hat{x},\epsilon}$, let $v_\lambda(\hat{x}) \in T_{\Phi_z^{-1}(\hat{x})} \wt{M}$ be given by
\begin{eqnarray}
\label{eq:v}v_\lambda(\hat{x}) &=& \left [ \frac{\partial}{\partial \hat{x}^j} \rho_\lambda \left (\Phi_z^{-1}(\hat{x}) \right ) \right ] \ g_z^{jk}(\Phi^{-1}_z(\hat{x})) \frac{\partial}{\partial \hat{x}^k} \\
&& \hskip.5in + \left (\sqrt{1 - \left |  \frac{\partial}{\partial \hat{x}^j} \rho_\lambda \left (\Phi_z^{-1}(\hat{x}) \right ) \right |_{g_z}^2 } \right ) \nu (\Phi_z^{-1}(\hat{x})) \nonumber
\end{eqnarray}
where $g_x^{jk}$ is the metric $g$ restricted to $\partial M$ expressed in the $\hat{x}$ coordinate frame. Then $\wt{\gamma}_{\hat{x},\epsilon}$ represents a segment of the geodesic $\gamma_{\Phi^{-1}_z(\hat{x}),-v_{\lambda_0}(\hat{x})}$. 

We now begin relating the constructions we have made so far with the data required to apply Theorem~\ref{thm2}. First, let us take a sufficiently small constant $\wt{\epsilon}>0$ and set
\[
x_0 = \gamma_{\wt{z},v_{\lambda_0}(0)}(\wt{\epsilon}), \quad \eta_0 = - \dot \gamma_{\wt{z},v_{\lambda_0}(0)}(\wt{\epsilon}).
\]
These can be constructed from the given data because $\gamma_{\wt{z},v_{\lambda_0}(0)}((0,\wt{\epsilon}) )\subset W' \cap (\wt{M} \setminus M)$ for $\wt{\epsilon}$ sufficiently small. We can also construct the map $\Phi^{-1}: U \rightarrow \wt{M}$ defined by
\[
\Phi^{-1}(\hat{x}) = \gamma_{\Phi^{-1}_z(\hat{x}),v_{\lambda_0}(\hat{x})}\left (\wt{\epsilon} + \rho_{\lambda_0}(\wt{z}) - \rho_{\lambda_0}(\Phi_z^{-1}(\hat{x})) \right ).
\]
Taking $\wt{\epsilon}$ sufficiently small, and possibly shrinking $U$ we can make $\Phi^{-1}$ a diffeomorphism onto its image, and so $\Phi$ is a coordinate map on the image of $\Phi^{-1}$ which, using the notation introduced earlier in the paper, we take to be $\Sigma_{0,t_0}$ with $t_0 := \wt{\epsilon}+ \rho_{\lambda_0}(\wt{z})$. Also following the earlier notation we let $\nu_{0,t_0}(\hat{x})$ be the outward pointing unit normal vector for $\Sigma_{0,t_0}$ at $\Phi^{-1}(\hat{x})$. 

Suppose now that $\hat{x} \in U$, and $t \in [\wt{\epsilon} + \rho_{\lambda_0}(\wt{z}) - \rho_{\lambda_0}(\Phi_z^{-1}(\hat{x})),t_0]$. If $\wt{\epsilon} + \rho_{\lambda_0}(\wt{z}) - \rho_{\lambda_0}(\Phi_z^{-1}(\hat{x}))$ and $t$ are not conjugate along $\gamma_{\Phi^{-1}(\hat{x}),-\nu_{0,t_0}(\hat{x})}$ then there exists a $\lambda \in \wt{\gamma}_{\hat{x},\epsilon}$ such that $\rho_\lambda(\Phi_z^{-1}(\hat{x})) = t- \wt{\epsilon} - \rho_{\lambda_0}(\wt{z})+ \rho_{\lambda_0}(\Phi_z^{-1}(\hat{x}))$. Now we have that $v_\lambda(\cdot)$ is defined on a neighborhood $U'$ of $\hat{x}$, and the set $\Sigma_{0,t}$ is then given by
\[
\Sigma_{0,t} = \left \{ \gamma_{\Phi_z^{-1}(\hat{x}'),v_\lambda(\hat{x}')} \left ( t- \rho_\lambda(\Phi_z^{-1}(\hat{x})) \right ) \ : \ \hat{x}' \in U' \right \}.
\]
Since we can recover this set and we know $g \in \wt{M}\setminus M$ we can calculate the shape operator $S^{t_0}_{0,t}(\hat{x})$ provided it exists. Therefore we may calculate $\mathbf{s}(\hat{x},t_0;0,t)$. Since there are only a finite number of $t \in [\wt{\epsilon} + \rho_{\lambda_0}(\wt{z}) - \rho_{\lambda_0}(\Phi_z^{-1}(\hat{x})),t_0]$ that are conjugate to $\wt{\epsilon} + \rho_{\lambda_0}(\wt{z}) - \rho_{\lambda_0}(\Phi_z^{-1}(\hat{x}))$ along $\gamma_{\Phi^{-1}(\hat{x}),-\nu_{0,t_0}(\hat{x})}$ we may in fact obtain $\mathbf{s}(\hat{x},t_0;0,t)$ for all $t \in [\wt{\epsilon} + \rho_{\lambda_0}(\wt{z}) - \rho_{\lambda_0}(\Phi_z^{-1}(\hat{x})),t_0]$ for which it is defined by continuity. Finally since $\gamma_{\Phi^{-1}(\hat{x}),-\nu_{0,t_0}(\hat{x})}((0,\wt{\epsilon} + \rho_{\lambda_0}(\wt{z}) - \rho_{\lambda_0}(\Phi_z^{-1}(\hat{x})))) \subset \wt{M} \setminus M$ we can also find $\mathbf{s}(\hat{x},t_0;0,t)$ for all $t \in (0,t_0]$ for which it is defined.

Now suppose that $t_0' \in \left (-\epsilon, \rho_{\lambda_0}(\Phi_z^{-1}(\hat{x}))\right )$. Provided that $t_0'$ and $\wt{\epsilon}$ are not conjugate along $\gamma_{x_0,\eta_0}$, there exists $\lambda_0' \in \wt{\gamma}_{\hat{x},\epsilon}$ such that $\wt{z} \in \mathrm{domain}(\rho_{\lambda_0'})$. We may now repeat the previous construction with $t_0$ replaced by $t_0'$ to obtain $\mathbf{s}(\hat{x},t_0';0,t)$ for all $t \in (0,t_0']$ where it is defined, and this gives all of the data necessary to apply Theorem~\ref{thm1}. Therefore we can recover the wave speed $v$ in a neighborhood of $\gamma_{\Phi^{-1}_z(\hat{x}),-v_{\lambda_0}(\hat{x})}\left ((0, t_0+\epsilon) \right )$. By the hypotheses $W'$ can be covered by sets of this form, and so we can recover the wave speed $v$ on all of $W'$ as claimed.
\endproof

\section{Conclusion}\label{sec:conc}

We generalized the method of Dix for reconstructing a depth varying
velocity in a half space, where depth is the Cartesian coordinate
normal to the boundary, to a procedure for reconstructing a conformally Euclidean metric on a region of $\mathbb{R}^n$ from expansions of diffraction travel times
generated by scatterers in the region and measured on its
boundary. Our procedure consists of two steps: In the first step, we
reconstruct the directional curvature operator along geodesics as well as
the metric in Riemannian normal coordinates. Riemannian normal
coordinates can be thought of as ``time'' coordinates as they appear
in so-called seismic time migration. We note that the directional curvature operator did not appear in the method of Dix because of the class of
velocity models he considered. In the second step, the velocity and
the geodesics on which the velocity is reconstructed are obtained
through a transformation to Cartesian coordinates; this can be thought of
as a generalization of the ``time-to-depth'' conversion in the framework
of Dix' original formulation. In dimension three or more both steps are essentially formulated in
terms of solving a closed system of nonlinear ordinary differential
equations, for example, by application of a Runge-Kutta
method. In dimension two the second step requires the solution of a Cauchy problem for an elliptic operator which may suffer from stability issues. Through the associated discretization, we accommodate the case
of a finite number of scatterers in the manifold. We admit the
formation of caustics.

\appendix

\section{Converting travel times measured at the boundary to the shape operator} \label{appendA}

In this appendix we show more explicitly how the somewhat abstract procedure described in the proof of Theorem~\ref{thm2} to move from travel times, or generalized distance functions, to the shape operators for wavefronts can be done in a particular case. We assume that $M$ is a lower half space and $v$ is constant in a neighborhood of the boundary and thus can be extended smoothly as a constant. Here we have the application in mind, and so we use the term travel time for the distance in the metric $g$.

\begin{figure}
\includegraphics[width=\textwidth]{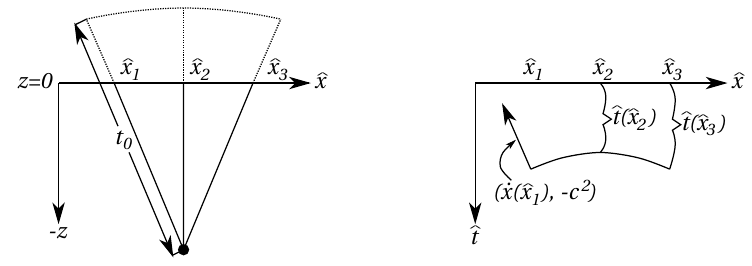}
\caption{On the left is illustrated a few rays originating at a common diffraction point and contributing to a single wavefront. The travel times are measured at $z=0$ producing a curve in $(\hat{x},\hat{t})$ space as illustrated on the right.}
\end{figure}

We use Cartesian coordinates $(x,z) \in \mathbb{R}^{n-1} \times \mathbb{R}$ and assume that we have a flat boundary at $\partial M = \{z=0\}$ so that $M = \{z \leq 0\}$ and $\wt{M} = \mathbb{R}^n$. We focus on a single diffraction point, $y_0$, and measure the travel times from $y_0$ as a function of the transverse coordinates at $\{z=0\}$ which we label as $\hat{x}$. We assume that the wave speed $v$ is equal to a constant $c$ in a neighborhood of $\{z = 0\}$ and write the travel time as $\hat{t}$ which is a function of $\hat{x}$ for $\hat{x}$ in some open set $U \subset \mathbb{R}^{n-1}$.

First we can determine the tangent vector of the ray passing through a certain point $\hat{x}$ of the boundary and corresponding to the travel time $\hat{t}(\hat{x})$ by
\[
\dot{x}^{k'}(\hat{x}) := \frac{\mathrm{d} x^{k'}}{\mathrm{d} t}(\hat{t}(\hat{x}),\hat{x}) = c^2 \frac{\partial \hat{t}}{\partial \hat{x}^{j'}}(\hat{x})\delta^{j'k'},
\]
\[
\dot{z}(\hat{x}) := \frac{\mathrm{d} z}{\mathrm{d} t}(\hat{t}(\hat{x}),\hat{x}) = \sqrt{c^2 - \left |\dot{x}(\hat{x}) \right |^2}
\]
(these equations should be compared with (\ref{eq:v})). We are using the prime notation here to indicate that the indices only run from $1$ to $n-1$. In practice it would be possible to calculate $\dot{x}$ by finding the normal to the surface $(\hat{x},\hat{t}(\hat{x}))$ which is proportional to the vector
\[
\left (-\dot{x}(\hat{x}) ,c^2\right ).
\]
Also note that if we look at multiple wavefronts coming from different diffraction points, then we can identify contributions from diffraction points along the same ray by the fact that they must satisfy
\[
\dot{x}_1(\hat{x}) = \dot{x}_2(\hat{x})
\]
(compare with (\ref{eq:wtgamma})) where $\dot{x}_1$ and $\dot{x}_2$ correspond to the two different wavefronts.

Now we choose $t_0 > \sup{\hat{t}}$ and consider the surface of points with travel time $t_0$ from the given diffraction point. In the body of the paper this surface is labeled $\Sigma_{0,t_0}$. We may parametrize the surface $\Sigma_{0,t_0}$ by $\hat{x}$ using the map (labeled $\Phi_{t_0}^{-1}$ in the body of the paper)
\[
\hat{x} \mapsto \left ( x = (t_0 - \hat{t}(\hat{x})) \dot{x}(\hat{x}) + \hat{x}, z = (t_0 - \hat{t}(\hat{x}) )\dot{z}(\hat{x})\right ).
\]
It is through this map that we define the coordinates on $\Sigma_{0,t_0}$. For $0\leq r \leq t_0-\hat{t}(\hat{x})$ we can also write a formula for the inverse of the coordinate map for the coordinates $(\hat{x},r)$ ($\Psi_{t_0}$ in the body of the paper). This is
\[
(\hat{x},r) \mapsto \left ( x = (t_0 - \hat{t}(\hat{x})-r) \dot{x}(\hat{x}) + \hat{x}, z = (t_0 - \hat{t}(\hat{x}) - r )\dot{z}(\hat{x}) \right ).
\]
Now we may calculate the matrices $F^k_j$ and $\widetilde{J}^k_j$ which give respectively the parallel translated fields $F_j$ and the coordinate vector fields $\partial/\partial \hat{x}^j$ (which are also Jacobi fields) in the Cartesian coordinates (note that $\widetilde{J}_j^k$ differs from $J_j^k$ in the body of the paper which represents the Jacobi fields with respect to the parallel frame).  To take advantage of the summation notation we also write $\hat{x}^n = r$ and $x^n = z$. Then we have
\[
\left . \frac{\partial}{\partial \hat{x}^j}\right |_{(\hat{x},r)} = \underbrace{\frac{\partial x^k}{\partial \hat{x}^j }(\hat{x},r)}_{\widetilde{J}_j^k(\hat{x},r,t_0)} \left .\frac{\partial}{\partial x^k} \right |_{(x,z)}
\]
where we may calculate (using the notation $s(\hat{x},r) = t_0 - \hat{t}(\hat{x}) - r$)
\[
\widetilde{J}(\hat{x},r,t_0) = \left (
\begin{matrix}
\delta^{k'}_{j'} +s \frac{\partial \dot{x}^{k'}}{\partial \hat{x}^{j'}} - \frac{1}{c^2}\dot{x}^{k'}\dot{x}^{l'} \delta_{l'j'} &
-\left \langle\dot{x}, \frac{\partial \dot{x}}{\partial \hat{x}^j} \right \rangle \frac{s}{\dot{z}} - \dot{x}^{l'} \delta_{l'j'} \frac{\dot{z}}{c^2}\\
-\dot{x}^{k'}& - \dot{z}
\end{matrix}
\right ).
\]
The parallel fields are the same as the coordinate vectors at $r =0$, and since for $0\leq r \leq t_0-\hat{t}(\hat{x})$ the wave speed is constant they are in fact constant in the Cartesian frame on the same set. Thus we have, for $0\leq r \leq t_0-\hat{t}(\hat{x})$,
\small
\[
F(\hat{x},r,t_0) = \left (
\begin{matrix}
\delta^{k'}_{j'} + \left (t_0 - \hat{t} \right ) \frac{\partial \dot{x}^{k'}}{\partial \hat{x}^{j'}} - \frac{1}{c^2}\dot{x}^{k'} \dot{x}^{l'}\delta_{l'j'} &
-\left \langle\dot{x}, \frac{\partial \dot{x}}{\partial \hat{x}^{j'}} \right \rangle \frac{\left (t_0-\hat{t} \right )}{\dot{z}} - \dot{x}^{l'} \delta_{l'j'} \frac{\dot{z}}{c^2}\\
-\dot{x}^{k'} & - \dot{z}
\end{matrix}
\right ).
\]
\normalsize
These matrices may also be written directly in terms of the travel time function $\hat{t}$. Indeed if we write $\mathrm{Hess}(\hat{t})$ for the Hessian matrix of $\hat{t}$, and $\nabla \hat{t}$ for the gradient then
\scriptsize
\[
\widetilde{J}(\hat{x},r,t_0) = \left (
\begin{matrix}
\delta +s \, c^2\mathrm{Hess}(\hat{t})- c^2 (\nabla\, \hat{t}) (\nabla \,\hat{t})^T &
-\frac{s \,c^2}{\sqrt{c^{-2}- \left |\nabla \hat{t} \right |^2}}  \mathrm{Hess}(\hat{t}) \nabla \, \hat{t} - c^2\sqrt{c^{-2}- \left |\nabla \hat{t} \right |^2} \nabla \, \hat{t}\\
-c^2 (\nabla \, \hat{t})^T& - c^2 \sqrt{c^{-2}- \left |\nabla \hat{t} \right |^2}
\end{matrix}
\right )
\]
\normalsize
and
\tiny
\[
F(\hat{x},r,t_0) = \left (
\begin{matrix}
\delta +(t_0 - \hat{t}) \, c^2\mathrm{Hess}(\hat{t})- c^2 (\nabla\, \hat{t}) (\nabla \,\hat{t})^T &
-\frac{(t_0 - \hat{t}) \,c^2}{\sqrt{c^{-2}- \left |\nabla \hat{t} \right |^2}}  \mathrm{Hess}(\hat{t}) \nabla \, \hat{t} - c^2\sqrt{c^{-2}- \left |\nabla \hat{t} \right |^2} \nabla \, \hat{t}\\
-c^2 (\nabla \, \hat{t})^T& - c^2 \sqrt{c^{-2}- \left |\nabla \hat{t} \right |^2}
\end{matrix}
\right ).
\]
\normalsize
Now we calculate the matrix $\mathbf{s}(\hat{x},t_0;,r,t_0)$ which corresponds with the shape operator of $\Sigma_{r,t_0}$ with respect to the parallel frame. To begin, recall that
\[
\mathbf{s}^k_j(\hat{x},t_0;,r,t_0) = \langle f^k(\hat{x},r), S_{r,t_0}^{t_0}(F_j(\hat{x},r) \rangle
\]
where $\{f^k\}_{k=1}^n$ is the dual frame for $\{F_j\}_{j=1}^n$. Now we can use the fact that the Christoffel symbols in Cartesian coordinates are zero where the wave speed is constant to calculate
\[
\begin{array}{lcl}
S_{r,t_0}^{t_0}(F_j(\hat{x},r)) & = & \nabla_{F_j} \left ( \dot{x}^{p'} \frac{\partial}{\partial x^{p'}} + \dot{z} \frac{\partial}{\partial z}\right )\\
& = & F_j^q\, \nabla_{\partial/\partial x^q} \left ( \dot{x}^{p'} \frac{\partial}{\partial x^{p'}} + \dot{z} \frac{\partial}{\partial z}\right )\\
& = & F_j^q \left ( \frac{\partial \dot{x}^{p'}}{ \partial x^q} \frac{\partial}{\partial x^{p'}} + \frac{\partial \dot{z}}{\partial x^q} \frac{\partial}{\partial z} \right )\\
& = & F_j^q \left ( \widetilde{J}^{-1}\right )^l_q \left ( \frac{\partial \dot{x}^{p'}}{ \partial \hat{x}^l} \frac{\partial}{\partial x^{p'}} + \frac{\partial \dot{z}}{\partial \hat{x}^l} \frac{\partial}{\partial z} \right )\\
& =  & F_j^q \left ( \widetilde{J}^{-1}\right )^l_q \left ( \frac{\partial \dot{x}^{p'}}{ \partial \hat{x}^l} \frac{\partial}{\partial x^{p'}} - \frac{1}{\dot{z}(\hat{x})}\left \langle\dot{x}(\hat{x}), \frac{\partial \dot{x}}{\partial \hat{x}^l}  \right \rangle \frac{\partial}{\partial z} \right )\\
& =  & F_j^q \left ( \widetilde{J}^{-1}\right )^l_q  \left (
\begin{matrix}  \frac{\partial \dot{x}}{ \partial \hat{x}} &
 - \frac{1}{\dot{z}} \left \langle\dot{x}, \frac{\partial \dot{x}}{\partial \hat{x}}  \right \rangle \\
 0 & 0
 \end{matrix} \right )_l^p (F^{-1})^k_p\, F_k.
\end{array}
\]
Therefore
\[
\mathbf{s}^k_j(\hat{x},t_0;,r,t_0) = F_j^q \left ( \widetilde{J}^{-1}\right )^l_q  \left (
\begin{matrix} \frac{\partial \dot x}{\partial \hat{x}} & 
 - \frac{1}{\dot{z}}  \frac{\partial \dot{x}}{\partial \hat{x}} \dot{x}\\ 0 & 0
 \end{matrix} \right )_l^p (F^{-1})^k_p.
\]
Once again, this can also be changed to an expression in terms of $\hat{t}$ and its derivatives as
\[
\mathbf{s}^k_j(\hat{x},t_0;r,t_0) = F \widetilde{J}^{-1}  \left (
\begin{matrix} c^2 \mathrm{Hess}(\hat{t}) &
-\frac{c^2}{\sqrt{c^{-2}- \left |\nabla \hat{t} \right |^2}} \mathrm{Hess}(\hat{t}) \nabla\, \hat{t}  \\
0 & 0 
 \end{matrix} \right ) F^{-1}
\]
and we may use this formula to calculate $\mathbf{s}^k_j(\hat{x},t_0;t_0-\hat{t}(\hat{x}),t_0)$ which we invert to give the initial conditions for the system \re{EQ A}. To calculate the derivatives with respect to $t$ we must consider multiple wavefronts, and using the comment above identify singularities in the wavefront corresponding to diffraction points along the same ray.

\end{document}